\documentclass[reqno, 10pt, oneside]{amsart}

\usepackage{latexsym, amsmath,amssymb}
\usepackage[hidelinks]{hyperref}
\ifdefined\pdfinterwordspaceon
  \ifnum\pdfoutput>0
    \input{glyphtounicode}
    \pdfinterwordspaceon
  \fi
\fi

\usepackage[initials, alphabetic,backrefs,msc-links]{amsrefs}
\usepackage{graphicx}
\usepackage{xcolor}
\usepackage{enumitem}

 \numberwithin{equation}{section}

\usepackage{amsmath}

\def\XXint#1#2#3{{\setbox0=\hbox{$#1{#2#3}{%
\int}$ }
\vcenter{\hbox{$#2#3$ }}\kern-.6\wd0}}

\renewcommand{\epsilon}{\varepsilon}
\newtheorem{theorem}{Theorem}

\newtheorem{lemma}[theorem]{Lemma}
\newtheorem{corollary}[theorem]{Corollary}

\newtheorem{proposition}[theorem]{Proposition}

\newtheorem{remark}[theorem]{Remark}
\newcommand{\bth}{\begin{theorem}}
\newcommand{\ble}{\begin{lemma}}
\newcommand{\bcor}{\begin{corr}}

\newcommand{\bdeff}{\begin{deff}}

\newcommand{\bprop}{\begin{proposition}}
\newcommand{\ele}{\end{lemma}}
\newcommand{\ecor}{\end{corr}}
\newcommand{\edeff}{\end{deff}}

\numberwithin{theorem}{section}

\newcommand{\eprop}{\end{proposition}}

\renewcommand{\Pi}{\varPi}

\renewcommand{\epsilon}{\varepsilon}

\newcommand{\R}{{\mathbb R}}

\DeclareMathOperator*{\essinf}{ess\,inf}

\newcommand{\norm}[1]{\left\|#1\right\|}

\usepackage{esint}

\usepackage{comment}

\allowdisplaybreaks

\def\vint_#1{\mathchoice%
        {\mathop{\kern 0.2em\vrule width 0.6em height 0.69678ex depth -0.58065ex
                \kern -0.8em \intop}\nolimits_{\kern -0.4em#1}}%
        {\mathop{\kern 0.1em\vrule width 0.5em height 0.69678ex depth -0.60387ex
                \kern -0.6em \intop}\nolimits_{#1}}%
        {\mathop{\kern 0.1em\vrule width 0.5em height 0.69678ex
            depth -0.60387ex
                \kern -0.6em \intop}\nolimits_{#1}}%
        {\mathop{\kern 0.1em\vrule width 0.5em height 0.69678ex depth -0.60387ex
                \kern -0.6em \intop}\nolimits_{#1}}}
\def\vintslides_#1{\mathchoice%
        {\mathop{\kern 0.1em\vrule width 0.5em height 0.697ex depth -0.581ex
                \kern -0.6em \intop}\nolimits_{\kern -0.4em#1}}%
        {\mathop{\kern 0.1em\vrule width 0.3em height 0.697ex depth -0.604ex
                \kern -0.4em \intop}\nolimits_{#1}}%
        {\mathop{\kern 0.1em\vrule width 0.3em height 0.697ex depth -0.604ex
                \kern -0.4em \intop}\nolimits_{#1}}%
        {\mathop{\kern 0.1em\vrule width 0.3em height 0.697ex depth -0.604ex
                \kern -0.4em \intop}\nolimits_{#1}}}

\newcommand{\aveint}[2]{\mathchoice%
        {\mathop{\kern 0.2em\vrule width 0.6em height 0.69678ex depth -0.58065ex
                \kern -0.8em \intop}\nolimits_{\kern -0.45em#1}^{#2}}%
        {\mathop{\kern 0.1em\vrule width 0.5em height 0.69678ex depth -0.60387ex
                \kern -0.6em \intop}\nolimits_{#1}^{#2}}%
        {\mathop{\kern 0.1em\vrule width 0.5em height 0.69678ex depth -0.60387ex
                \kern -0.6em \intop}\nolimits_{#1}^{#2}}%
        {\mathop{\kern 0.1em\vrule width 0.5em height 0.69678ex depth -0.60387ex
                \kern -0.6em \intop}\nolimits_{#1}^{#2}}}

\def\XXint#1#2#3{{\setbox0=\hbox{$#1{#2#3}{\int}$}
    \vcenter{\hbox{$#2#3$}}\kern-.5\wd0}}

\newcommand{\vertiii}[1]{{\left\vert\kern-0.25ex\left\vert\kern-0.25ex\left\vert #1 
    \right\vert\kern-0.25ex\right\vert\kern-0.25ex\right\vert}}

\newcommand{\vertii}[1]{{\left\vert\kern-0.25ex\left\vert\kern-0.25ex  #1 
    \kern-0.25ex\right\vert\kern-0.25ex\right\vert}}

\makeatletter
\@namedef{subjclassname@2020}{\textup{2020} Mathematics Subject Classification}
\makeatother

\begin{document}

\title[Quantitative BMO--BLO Estimates for the Maximal Function]
{Quantitative BMO--BLO Estimates for the Hardy--Littlewood Maximal Function
}

\author[A. Claros]{Alejandro Claros}

\address[A. Claros]{BCAM -- Basque Center for Applied Mathematics, Bilbao, Spain\newline
Universidad del Pa\'{\i}s Vasco / Euskal Herriko Unibertsitatea (UPV/EHU), Bilbao, Spain}

\email{aclaros@bcamath.org, aclaros003@ikasle.ehu.eus}

\thanks{This research is supported by the Basque Government through the BERC 2022--2025 program, by the Ministry of Science and Innovation through Grant PRE2021-099091 funded by BCAM Severo Ochoa accreditation CEX2021-001142-S/MICIU/AEI/10.13039/501100011033 and by ESF+, and by the project PID2023-146646NB-I00 funded by MICIU/AEI/10.13039/501100011033 and by ESF+.}

\subjclass[2020]{Primary 42B25; Secondary 42B35}



\keywords{BMO, John--Nirenberg inequality, Muckenhoupt weights, maximal function}

\begin{abstract}
In this note, we study a quantitative extension of the John--Nirenberg inequality for the Hardy--Littlewood maximal function of a $\operatorname{BMO}$ function. More precisely, for every nonconstant locally integrable function $f$ such that $Mf$ is not identically infinite, we prove the inequality
\begin{equation*}
	\left( \frac{1}{w(Q)}\int_Q \left( \frac{Mf(x) - \essinf_{Q} Mf }{M^\# f(x)} \right)^p w(x)\,dx\right)^\frac{1}{p} \le c_n \, [w]_{A_\infty}\,  p
\end{equation*}
for every cube $Q$, every $1\le p<\infty$ and every weight $w\in A_\infty$, where $[w]_{A_\infty}$ denotes the Fujii--Wilson $A_\infty$ constant. This result extends the classical boundedness $\|Mf\|_{\operatorname{BLO}}\le C_n\|f\|_{\operatorname{BMO}}$ (\cites{BDS, Bennett}). Furthermore, we show that the class $A_\infty$ is both necessary and sufficient for this inequality to hold, providing a new characterization of $A_\infty$ in terms of the action of the maximal operator on bounded oscillation spaces.
\end{abstract}

\maketitle

\section{Introduction and statement of the main results}

The theory of functions of bounded mean oscillation ($\operatorname{BMO}$) plays a fundamental role in modern harmonic analysis. The classical John--Nirenberg theorem \cite{JN} states that functions of bounded mean oscillation are locally exponentially integrable. More precisely, if $f \in \operatorname{BMO}$, then there exist dimensional constants $c_1,c_2>0$ such that for every cube $Q\subset \R^n$,
\begin{equation*}
	|\{x\in Q:\, |f(x)-f_Q|>t\}|\le c_1 e^{-c_2 t/\|f\|_{\operatorname{BMO}}}|Q|,\qquad t>0,
\end{equation*}
where
\begin{equation*}
	\|f\|_{\operatorname{BMO}}=\sup_Q \frac{1}{|Q|}\int_Q |f(x)-f_Q| \,dx
\end{equation*}
with $f_Q=\frac{1}{|Q|}\int_Q f$, and the supremum is taken over all cubes in $\R^n$ with sides parallel to the coordinate axes. This inequality establishes that bounded mean oscillation implies a form of exponential control stronger than any $L^p$ bound. This property also identifies $\operatorname{BMO}$ as a natural endpoint substitute for $L^\infty$ in the scale of Lebesgue spaces and as a central object in the study of singular integrals, interpolation, and partial differential equations.

The Hardy--Littlewood maximal operator
\begin{equation*}
	Mf(x)=\sup_{Q\ni x}\frac{1}{|Q|}\int_Q |f(y)|\,dy,
\end{equation*}
and the sharp maximal operator introduced by Fefferman and Stein \cite{FS},
\begin{equation*}
	M^\#f(x)=\sup_{Q\ni x}\frac{1}{|Q|}\int_Q |f(y)-f_Q|\,dy,
\end{equation*}
are key tools in the study of singular integrals, regularity, weighted inequalities, and self-improvement phenomena.  These operators satisfy the pointwise inequality $M^\#f(x)\le 2 Mf(x)$, and the Fefferman--Stein inequality shows that both operators are equivalent in $L^p$ norms for $1<p<\infty$ and $f\in L^p$. The exponential integrability in the John--Nirenberg inequality may be viewed as an endpoint version of this relation.

Many authors have studied quantitative versions of the John--Nirenberg inequality, including extensions to weighted settings \cites{MRR, PR, OPRR, JCCP}, to metric measure spaces \cites{Stromberg, StrombergTorchinsky, FPW, DGL, Kim, GK}, and even to more general measures and capacities \cites{ChenSpector, BCRS}. Indeed, the weighted John--Nirenberg, Poincar\'e, and Poincar\'e--Sobolev inequalities can be unified within a single theory through the so-called self-improving phenomena (see \cites{FPW, MP, PR, GLP}).

In \cite{Karagulyan}, Karagulyan obtained an extension of the John--Nirenberg inequality, later refined by Canto and P\'erez \cite{JCCP}. This extension shows that one can recover the exponential decay of the oscillation of any function once it is normalized by its sharp maximal function. More precisely, the main result of \cite{JCCP} established the inequality
\begin{equation}\label{eq: extension JN}
	\left(\frac{1}{w(Q)}\int_Q \left(\frac{M_Q(f-f_Q)(x)}{M^\#f(x)}\right)^p
w(x)\,dx\right)^{\frac{1}{p}} \le C_n [w]_{A_\infty} p,
\end{equation}
for every cube $Q$, every $1\le p<\infty$ and every weight $w\in A_\infty$ (we also refer to \cite{LLO} for a different proof of \eqref{eq: extension JN}). Here, $M_Q$ is the local dyadic maximal operator, and $[w]_{A_\infty}$ denotes the Fujii--Wilson $A_\infty$ constant, introduced in \cite{HP}, defined by
\begin{equation*}
	[w]_{A_\infty}=\sup_Q \frac{1}{w(Q)}\int_Q M(w\chi_Q)(x)\,dx.
\end{equation*}

In a classical result, Bennett, DeVore, and Sharpley \cite{BDS} proved that the Hardy--Littlewood maximal operator $M$ maps $\operatorname{BMO}$ into itself. Shortly after, Bennett \cite{Bennett} refined this result by proving that $M$ actually maps $\operatorname{BMO}$ into a smaller space, the space of bounded lower oscillation functions ($\operatorname{BLO}$), satisfying
\begin{equation}\label{eq:M maps BMO into BLO}
	\|Mf\|_{\operatorname{BLO}} \le C_n \|f\|_{\operatorname{BMO}},
\end{equation}
for every function $f\in \operatorname{BMO}$ such that $Mf$ is not identically infinite. Introduced by Coifman and Rochberg in \cite{CR}, the space $\operatorname{BLO}$ consists of all locally integrable functions such that
\begin{equation*}
	\|f\|_{\operatorname{BLO}}=\sup_Q \frac{1}{|Q|}\int_Q \left(f(x)-\essinf_{Q}f\right)\,dx<\infty.
\end{equation*}

From the perspective of weight theory, these spaces admit a natural characterization in terms of logarithms of Muckenhoupt weights, as shown in \cite{CR}. In fact, one can write
\begin{equation*}
	\operatorname{BLO} = \{ \alpha \log w : w \in A_1,\, \alpha >0\},
\end{equation*}
and
\begin{equation*}
\operatorname{BMO} = \{ \alpha \log w : w \in A_\infty,\, \alpha \in \R \}.
\end{equation*}
This formulation shows the intrinsic link between bounded oscillation spaces and Muckenhoupt weights, and clearly shows the strict inclusion $\operatorname{BLO} \subsetneq \operatorname{BMO}$.

The present note is devoted to proving a quantitative extension of the mapping property \eqref{eq:M maps BMO into BLO}, analogous to \eqref{eq: extension JN}. Our proof combines the local estimate of Canto and P\'erez \cite{JCCP}*{Theorem 5.1} with Bennett's treatment of the nonlocal part \cite{Bennett}, and then passes to $A_\infty$ weights by means of the sharp reverse H\"older inequality \cite{HPR}. Concretely, we prove the following result.

\begin{theorem}\label{thm: 1}
	Let $f$ be a nonconstant locally integrable function such that $Mf$ is not identically infinite, and let $w\in A_\infty$ be a weight. Then, for any cube $Q$ and $1\le p <\infty$, we have
	\begin{equation}\label{eq: main eq}
		\left( \frac{1}{w(Q)}\int_Q \left( \frac{Mf(x) - \essinf_{Q} Mf }{M^\# f(x)} \right)^p w(x)\, dx\right)^\frac{1}{p} \le c_n\,  [w]_{A_\infty}\, p.
	\end{equation}
	In particular, if $f\in \operatorname{BMO}$, then for each $1\le p<\infty$ we have
	\begin{equation}\label{eq: M map BMO into BLOw}
		\left( \frac{1}{w(Q)}\int_Q \left( Mf(x) - \essinf_{Q} Mf  \right)^p w(x)\, dx\right)^\frac{1}{p} \le c_n \, [w]_{A_\infty} \, p \,  \|f\|_{\operatorname{BMO}},
	\end{equation}
	for all cubes $Q$. 
\end{theorem}

\begin{remark}
The assumption that $Mf$ is not identically infinite cannot be removed. In \cite{BDS}, it is shown that this condition is necessary for the validity of the boundedness of $M$ on $\operatorname{BMO}$, and an explicit counterexample is provided by the function $f(x)=\log|x|$. We also refer to the recent work \cite{BGSKM}*{Theorem A.1}, which characterizes the class of measures $\mu$ for which $M\mu$ is finite almost everywhere (see also \cite{FK00} for the corresponding class of locally integrable functions).
\end{remark}

The linear dependence on $p$ in the right-hand side of \eqref{eq: main eq} yields exponential estimates as a direct corollary. 

\begin{corollary}\label{cor: 1}
	Let $f$ be a nonconstant locally integrable function such that $Mf$ is not identically infinite. Let $w\in A_\infty$ be a weight. 
	\begin{itemize}
		\item There exist dimensional constants $c_1,c_2>0$ such that 
	\begin{equation}\label{eq: weak JN maximal}
		w \left( \left\lbrace x\in Q : \frac{Mf(x) - \essinf_{Q} Mf }{M^\# f(x)} > t \right \rbrace \right) \le c_1 e^{-\frac{c_2 t}{[w]_{A_\infty}}} w(Q),
	\end{equation}
	for all $t>0$ and all cubes $Q$. That is, 
	\begin{equation*}
		\norm{\frac{Mf - \essinf_{Q} Mf }{M^\# f}}_{\exp L \left(Q , \frac{w(x)\, dx}{w(Q)}\right)} \le c_n [w]_{A_\infty}. 
	\end{equation*}
	\item There exist dimensional constants $c_1,c_2>0$ such that for every cube $Q$ and for every $\lambda, \gamma>0$ the following good-lambda inequality holds:
	\begin{equation}\label{eq: good lambda}
		w\left( \left\lbrace x\in Q : Mf(x) - \essinf_{Q} Mf > \lambda, M^\# f(x) \le \gamma \lambda \right \rbrace \right) \le c_1 e^{-\frac{c_2 }{ [w]_{A_\infty} \gamma}} w(Q).
	\end{equation}
	\end{itemize}
\end{corollary}

\begin{remark}
	Inequality \eqref{eq: good lambda} improves, in its dependence on $\gamma$, the good-lambda estimate of Journ\'e \cite{J}*{p. 43}, which (for the dyadic versions of $M$ and $M^\#$ and a dyadic cube $Q$) gives
	\begin{align*}
		|\{x\in Q :{} & Mf(x)-\essinf\nolimits_Q Mf>2\lambda,\ M^\# f(x)<\gamma\lambda\}| \\
		&\qquad    \le 2^n\gamma\,|\{x\in Q : Mf(x)-\essinf\nolimits_Q Mf>\lambda\}| \le 2^n\gamma\,|Q|,
	\end{align*}
	with linear decay in $\gamma$. In the weighted case, combining this with the $A_\infty$ property recalled in Section 3 only yields a bound of order $\gamma^{1/(2^{n+1}[w]_{A_\infty})}\,w(Q)$, whereas \eqref{eq: good lambda} gives $c_1e^{-c_2/([w]_{A_\infty}\gamma)}\,w(Q)$.
\end{remark}

The previous inequality \eqref{eq: M map BMO into BLOw}, in the particular case $p=1$, tells us that if $w\in A_\infty$ then we have the boundedness
\begin{equation*}
	\|Mf\|_{\operatorname{BLO}_w} \le C_n [w]_{A_\infty} \|f\|_{\operatorname{BMO}},
\end{equation*}
where 
\begin{equation*}
	\|f\|_{\operatorname{BLO}_w} := \sup_Q \frac{1}{w(Q)}\int_Q \left( f(x) - \essinf_{Q} f  \right) w(x)\, dx.
\end{equation*}
We can also prove that the converse implication holds: if the above inequality is satisfied, with some finite constant in place of $C_n[w]_{A_\infty}$, for all locally integrable functions $f$ for which $Mf$ is not identically infinite, then the weight $w$ necessarily belongs to the class $A_\infty$. This equivalence yields a new quantitative characterization of the $A_\infty$ class in terms of the boundedness of the Hardy--Littlewood maximal operator acting on $\operatorname{BMO}$.

\begin{theorem}\label{thm: characterization Ainfty}
	Let $w$ be a weight. Then there exist some dimensional constants $c_n, C_n$ such that 
	\begin{equation}\label{eq: equiv}
		c_n [w]_{A_\infty} \le \sup_{f} \| Mf\|_{\operatorname{BLO}_w } \le C_n [w]_{A_\infty}, 
	\end{equation}
	where the supremum is taken over all functions $f$ such that $\|f\|_{\operatorname{BMO}}=1$ and $Mf$ is not identically infinite. As a consequence, we have the following quantitative equivalence of the Fujii--Wilson $A_\infty$ constant, 
	\begin{equation*}
		[w]_{A_{\infty}} \simeq_n \sup_{f} \sup_Q \frac{1}{w(Q)}\int_Q \left( Mf(x) - \essinf_{Q} Mf \right) w(x)\, dx .
	\end{equation*}
\end{theorem}

As a direct corollary, we obtain that the class $A_\infty$ is both necessary and sufficient for the validity of inequality \eqref{eq: main eq}. 

\begin{corollary} \label{cor: char}
Let $w$ be a weight. 
\begin{itemize}
	\item Let $f$ be a nonconstant locally integrable function such that $Mf$ is not identically infinite. If $w\in A_\infty$, then for any cube $Q$ and $1\le p <\infty$, we have
	\begin{equation*}
		\left( \frac{1}{w(Q)}\int_Q \left( \frac{Mf(x) - \essinf_{Q} Mf }{M^\# f(x)} \right)^p w(x)\,dx\right)^\frac{1}{p} \le c_n [w]_{A_\infty} p.
	\end{equation*}
	\item Assume that the estimate 
	 \begin{equation}\label{eq: cor 2}
		\left( \frac{1}{w(Q)}\int_Q \left( \frac{Mf(x) - \essinf_{Q} Mf }{M^\# f(x)} \right)^p w(x)\, dx\right)^\frac{1}{p} \le C\,  p
	\end{equation}
	holds for every nonconstant locally integrable function $f$ such that $Mf$ is not identically infinite, every $1\le p<\infty$, and every cube $Q$. Then $w\in A_\infty$ and $[w]_{A_\infty}\lesssim_n C$. 
\end{itemize}
\end{corollary}

\begin{remark}
A similar characterization holds for the main result of \cite{JCCP}, namely inequality \eqref{eq: extension JN}. The proof follows exactly the same argument as in the corollary above, replacing the use of Theorem \ref{thm: characterization Ainfty} with \cite{OPRR}*{Theorem 1.2}.
\end{remark}

As usual, $C$ denotes a positive constant, possibly varying from line to line. We write $C_{\alpha, \beta, ...}$ to denote a constant depending only on $\alpha, \beta,...$.

\section{The unweighted case}

In this section, we present the proof of the main result, Theorem \ref{thm: 1}, in the unweighted setting. We first recall the following theorem from \cite{JCCP}*{Theorem 5.1}, which extends the exponential-type form of \eqref{eq: extension JN} to the classical (non-dyadic) maximal function. The proof relies on a result of \cite{Conde}, which is used to combine estimates over a finite family of different dyadic lattices and then recover the non-dyadic case.

\begin{theorem}
	Let $Q$ be an arbitrary cube and $f$ a nonconstant locally integrable function. Then for any $\lambda>0$ we have
\begin{equation}\label{eq: weak extJN}
	\left|\left\{x \in Q: \frac{M\left(\left(f-f_Q\right) \chi_Q\right)(x)}{M^\# f(x)}>\lambda\right\}\right| \leq C e^{-c \lambda}|Q|,
\end{equation}
where $C, c>0$ are dimensional constants.
\end{theorem}

Using the layer-cake formula, \eqref{eq: weak extJN}, and noting that $\Gamma(p+1)^{1/p} \lesssim p$, we obtain the following inequality,
\begin{equation}\label{eq: nondyadic extJN}
	\left(\frac{1}{|Q|}\int_Q \left(\frac{M((f-f_Q)\chi_Q)(x)}{M^\#f(x)}\right)^p \, dx\right)^{\frac{1}{p}} \le C_n  p,
\end{equation}
for each cube $Q$, and each $1\le p<\infty$ (we refer to the proof of Theorem \ref{thm: 1} below for more details). With this preliminary estimate established, we are now in a position to prove the unweighted version of Theorem \ref{thm: 1}.

\begin{proposition}\label{Thm1}
	Let $f$ be a nonconstant locally integrable function such that $Mf$ is not identically infinite. Then, there exists a dimensional constant $c_n>0$ such that for every $1\le p <\infty$, we have
	\begin{equation*}
		\left( \frac{1}{|Q|}\int_Q \left( \frac{Mf(x) - \essinf_{Q} Mf }{M^\# f(x)} \right)^p \, dx\right)^\frac{1}{p} \le c_n \, p,
	\end{equation*}
	for every cube $Q$.
\end{proposition}

In the proof, we decompose $Mf$ into two terms: the local and nonlocal parts. The local contribution is estimated using \eqref{eq: nondyadic extJN}, whereas the nonlocal term is treated following ideas from \cite{Bennett}.

\begin{proof}
	Fix a cube $Q$ and $1\le p<\infty$. We decompose the maximal function as
	\begin{align*}
		Mf(x) = & M( (f-f_{3Q})\chi_{3Q} + f_{3Q}\chi_{3Q} + f\chi_{(3Q)^c})(x)\\
		\le & M( (f-f_{3Q})\chi_{3Q})(x) + Mg(x),
	\end{align*}
	where $g(x):= f_{3Q}\chi_{3Q}(x) + f(x)\chi_{(3Q)^c}(x)$, and we denote by $3Q$ the cube with the same center as $Q$ and side length $3\ell(Q)$. Therefore,
	\begin{align*}
		Mf(x) - \essinf_{Q} Mf \le & M( (f-f_{3Q})\chi_{3Q})(x) + Mg(x) - \essinf_{Q} Mf\\
		\le & M( (f-f_{3Q})\chi_{3Q})(x) + \left( Mg(x) - \essinf_{Q} Mf \right)^+,
	\end{align*}
	where $a^+=\max (a, 0)$. This yields the estimate 
	\begin{align*}
		\left( \frac{1}{|Q|}\int_Q  \left( \frac{Mf - \essinf_{Q} Mf }{M^\# f} \right)^p \, dx\right)^\frac{1}{p} \le & \left( \frac{1}{|Q|}\int_Q \left( \frac{M( (f-f_{3Q})\chi_{3Q})}{M^\# f} \right)^p \, dx\right)^\frac{1}{p}\\
		&+ \left( \frac{1}{|Q|}\int_Q  \left( \frac{(Mg - \essinf_{Q} Mf)^+ }{M^\# f} \right)^p \, dx\right)^\frac{1}{p} \\
		=: & I_1 + I_2. 
	\end{align*}
	
	By applying \eqref{eq: nondyadic extJN} over the larger cube $3Q$, we obtain
	\begin{align*}
		I_1 \le & \left( \frac{1}{|Q|}\int_{3Q} \left( \frac{M( (f-f_{3Q})\chi_{3Q})(x)}{M^\# f(x)} \right)^p \, dx\right)^\frac{1}{p}\\
		\le & \left( \frac{3^n}{|3Q|}\int_{3Q} \left( \frac{M( (f-f_{3Q})\chi_{3Q})(x)}{M^\# f(x)} \right)^p\,  dx\right)^\frac{1}{p}\\
		\le & c_n  p.
	\end{align*}
	
	To estimate the second term, we bound $Mg$ pointwise in $Q$. Let $x\in Q$ and $P$ be a cube such that $x\in P$. We consider two cases. First, if $P\subseteq 3Q$, then
	\begin{align*}
		\frac{1}{|P|}\int_P |g(y)|\, dy = |f_{3Q}| \le \frac{1}{|3Q|}\int_{3Q} |f(y)|\, dy \le \essinf_{z\in Q} Mf(z). 
	\end{align*}
	Otherwise, if $P\cap (3Q)^c\neq \emptyset$, we consider the smallest cube $\bar{P}$ that contains both $P$ and $3Q$. Since $x\in Q\cap P$, we have $|\bar{P}|\le c_n |P|$. Hence,
	\begin{align*}
		\frac{1}{|P|}\int_P |g(y)|\, dy \le & \frac{1}{|P|}\int_P |g(y)-f_{\bar{P}}|\, dy + |f_{\bar{P}}| \\
		=: & II_1 + II_2. 
	\end{align*}
	We can estimate $II_2$ trivially,
	\begin{align*}
		II_2 \le \frac{1}{|\bar{P}|}\int_{\bar{P}} |f(y)|\, dy \le \essinf_{z\in Q} Mf(z). 
	\end{align*}
	Using the definition of $g$, we have
	\begin{align*}
		II_1 \le & c_n \frac{1}{|\bar{P}|}\int_{\bar{P}} |g(y)-f_{\bar{P}}|\, dy \\
		= & c_n \frac{1}{|\bar{P}|}\int_{\bar{P}} |f_{3Q}\chi_{3Q}(y) + f(y)\chi_{(3Q)^c}(y)-f_{\bar{P}}|\, dy \\
		= & c_n \frac{1}{|\bar{P}|}\int_{3Q} |f_{3Q}\chi_{3Q}(y) + f(y)\chi_{(3Q)^c}(y)-f_{\bar{P}}|\, dy \\
		& +c_n \frac{1}{|\bar{P}|}\int_{\bar{P}\cap (3Q)^c} |f_{3Q}\chi_{3Q}(y) + f(y)\chi_{(3Q)^c}(y)-f_{\bar{P}}|\, dy \\
		= & c_n \frac{1}{|\bar{P}|}\int_{3Q} |f_{3Q}-f_{\bar{P}}|\, dy  +c_n \frac{1}{|\bar{P}|}\int_{\bar{P}\cap (3Q)^c} |f(y)-f_{\bar{P}}|\, dy \\
		\le & c_n \frac{1}{|\bar{P}|} |3Q| \frac{1}{|3Q|} \int_{3Q} |f(y)-f_{\bar{P}}|\, dy + c_n \frac{1}{|\bar{P}|}\int_{\bar{P}\cap (3Q)^c} |f(y)-f_{\bar{P}}|\, dy\\
		= & c_n \frac{1}{|\bar{P}|}\int_{\bar{P}} |f(y)-f_{\bar{P}}|\, dy\\
		\le & c_n M^\# f(x). 
	\end{align*}
	
	We have proved 
	\begin{equation*}
		\frac{1}{|P|}\int_P |g(y)|\, dy \le c_n M^\# f(x) + \essinf_{z\in Q} Mf(z),
	\end{equation*}
	for each cube $P$ that contains $x$, and thus
	\begin{equation*}
		Mg(x) \le c_n M^\# f(x) + \essinf_{z\in Q} Mf(z),
	\end{equation*}
	for every $x\in Q$. Using this estimate, we can bound $I_2$, 
	\begin{align*}
		I_2 \le & \left( \frac{1}{|Q|}\int_Q  \left( \frac{c_n M^\# f(x)  }{M^\# f(x)} \right)^p \, dx\right)^\frac{1}{p} =c_n \le c_n p.
	\end{align*}
	This concludes the proof. 
\end{proof}

The linear dependence on $p$ in Proposition \ref{Thm1} implies exponential-type integrability (see \cite{JCCP}*{Lemma 2.3}).

\begin{corollary}\label{cor: weak unweighted}
	Let $f$ be a nonconstant locally integrable function such that $Mf$ is not identically infinite. Then, there exists a dimensional constant $c_n>0$ such that 
	\begin{equation}\label{eq: weak unweighted estimate}
		\left| \left\lbrace x\in Q : \frac{Mf(x) - \essinf_{Q} Mf }{M^\# f(x)} > t \right \rbrace \right| \le 2 e^{-c_n t} |Q|,
	\end{equation}
	for all $t>0$ and all cubes $Q$. 
\end{corollary}

\section{Weighted estimates and the characterization of \texorpdfstring{$A_\infty$}{A-infinity}}

A weight $w \in A_\infty$ satisfies
\begin{equation*}
    \frac{w(E)}{w(Q)} \le 2 \left( \frac{|E|}{|Q|} \right)^{\frac{1}{2^{n+1}[w]_{A_\infty}}},
\end{equation*}
for every cube $Q$ and every measurable set $E \subseteq Q$ (this follows from \cite{HPR}*{Theorem 2.3} and H\"older's inequality). Using this property of $A_\infty$ weights and \eqref{eq: weak unweighted estimate}, we can prove Corollary \ref{cor: 1}. We have chosen to follow the strategy of first establishing the unweighted case and then passing to the weighted estimate using the above property. This approach avoids relying on the doubling property of $A_\infty$ weights, whose doubling constant is known to depend doubly exponentially on $[w]_{A_\infty}$ (see \cite{HaPa}*{p. 945}). Our approach yields linear dependence on $[w]_{A_\infty}$ in \eqref{eq: main eq}.

\begin{proof}[Proof of Corollary \ref{cor: 1}]
	Let $w\in A_\infty$ and fix a cube $Q$. Let $t>0$ and denote 
	\begin{equation*}
		E_t = \left\lbrace x\in Q : \frac{Mf(x) - \essinf_{Q} Mf }{M^\# f(x)} > t \right \rbrace.
	\end{equation*}
	
	By the previous property of $A_\infty$ weights and \eqref{eq: weak unweighted estimate},
	\begin{align*}
		w \left( E_t \right) \le & 2 \left( \frac{|E_t|}{|Q|}\right)^\frac{1}{2^{n+1}[w]_{A_\infty}} w(Q)\\
		\le & 2 \left( 2 e^{-c_n t}\right)^\frac{1}{2^{n+1}[w]_{A_\infty}} w(Q)\\
		\le & 2^{1+\frac{1}{2^{n+1}}} e^{\frac{-c_n t}{2^{n+1}[w]_{A_\infty}}}w(Q).
	\end{align*}
	This concludes the proof of \eqref{eq: weak JN maximal}. To prove \eqref{eq: good lambda}, it is enough to observe that
	\begin{equation*}
		\left\lbrace x\in Q : Mf(x) - \essinf_{Q} Mf > \lambda, M^\# f(x) \le \gamma \lambda \right \rbrace \subseteq E_\frac{1}{\gamma},
	\end{equation*}
	and apply \eqref{eq: weak JN maximal}. 
\end{proof}

\begin{proof}[Proof of Theorem \ref{thm: 1}]
	Let $1\le p<\infty$, $w\in A_\infty$ and let $Q$ be a cube. For each $t>0$, define
	\begin{equation*}
		E_t = \left\lbrace x\in Q : \frac{Mf(x) - \essinf_{Q} Mf }{M^\# f(x)} > t \right \rbrace.
	\end{equation*}
	
	By the layer-cake formula and \eqref{eq: weak JN maximal}, we have
	\begin{align*}
		\left( \frac{1}{w(Q)}\int_Q \left( \frac{Mf(x) - \essinf_{Q} Mf }{M^\# f(x)} \right)^p w(x)\, dx\right)^\frac{1}{p}  = & \left( p \int_0^\infty t^{p-1} \frac{w(E_t)}{w(Q)} \, dt \right) ^\frac{1}{p}\\
		\le & \left( c_1 p \int_0^\infty t^{p-1} e^{-\frac{c_2 t}{[w]_{A_\infty}}} \, dt \right) ^\frac{1}{p}\\
		= & c_1^\frac{1}{p} \frac{[w]_{A_\infty}}{c_2} \left(  p \int_0^\infty s^{p-1} e^{-s} \, ds \right) ^\frac{1}{p}\\
		= &  c_1^\frac{1}{p} \frac{[w]_{A_\infty}}{c_2} \left(  p \Gamma(p) \right) ^\frac{1}{p}\\
		\le & \frac{c_1}{c_2} [w]_{A_\infty} \Gamma(p+1)^\frac{1}{p},
	\end{align*}
	 where we used that $c_1\ge 1$, and we obtain the desired estimate by noting that $\Gamma(p+1)^\frac{1}{p}\lesssim p$ (see for instance \cite{Grafakos}*{Corollary 6.2.5}).
	 \end{proof}

\begin{remark}
Using similar arguments, one can obtain Lorentz-norm variants of Theorem \ref{thm: 1}, but with a different dependence on $p$: for $1\le q<\infty$, the $L^{p,q}(Q,\frac{w\,dx}{w(Q)})$-norm of $\frac{Mf-\essinf_Q Mf}{M^\# f}$ is bounded by $c_{n,q}[w]_{A_\infty}\,p^{1+1/q}$. We do not pursue that direction.
\end{remark}

The remainder of this section is devoted to the proof of Theorem \ref{thm: characterization Ainfty}.  Within the proof, we will address some technical difficulties caused by the maximal operator, and then follow the lines of the proof of Theorem 1.2 in \cite{OPRR}.

\begin{proof}[Proof of Theorem \ref{thm: characterization Ainfty}]
	In \eqref{eq: M map BMO into BLOw} (with $p=1$) we have already proved that
	\begin{equation*}
		\|Mf\|_{\operatorname{BLO}_w } \le c_n [w]_{A_\infty}\|f\|_{\operatorname{BMO}}, 
	\end{equation*}
	which provides the second inequality in \eqref{eq: equiv}. To prove the first inequality, we define
	\begin{equation*}
		X = \sup_{\|f\|_{\operatorname{BMO}}=1} \| Mf\|_{\operatorname{BLO}_w } ,
	\end{equation*}
	where the supremum is taken over all functions $f$ such that $\|f\|_{\operatorname{BMO}}=1$ and $Mf$ is not identically infinite. We have to prove $[w]_{A_\infty} \lesssim_n X$. Without loss of generality, we can assume that $X<\infty$. 
	
	We first establish the inequality $X \ge c_n$. Let $Q_0= [0,1]^n$ be the unit cube and decompose it into its $2^n$ dyadic children, and let $E$ be one of them satisfying $w(E)=\max\{w(Q_j): Q_j \text{ is a dyadic child of }Q_0\}$. Then, by the additivity of the measure $w$, we have $w(E)\ge 2^{-n}w(Q_0)$.  Let $f:=2\chi_E$. By \cite{BMOnormcharacteristic}, we have $\|\chi_E\|_{\operatorname{BMO}} = \tfrac{1}{2}$, and thus $\|f\|_{\operatorname{BMO}}=1$. For every $x\in Q_0$ we have
	\begin{equation*}
		M f(x)\ge \frac1{|Q_0|}\int_{Q_0} f(y)\,dy = \frac{2|E|}{|Q_0|} = 2^{1-n}.
	\end{equation*}
	Therefore, we have proved $\essinf_{x\in Q_0} M f(x)\ge 2^{1-n}$. To prove the reverse inequality, let $z_0$ be the vertex of $Q_0$ strictly opposite to $E$. That is, the distance in the supremum norm between $z_0$ and $E$ is exactly $\tfrac{1}{2}$. Let $P$ be any cube of side length $L$ and such that $z_0\in P$. If $L \le \tfrac{1}{2}$, the interior of $P$ cannot intersect $E$, and the average of $f$ over $P$ is $0$. If $L > \tfrac{1}{2}$, the measure of the intersection is bounded by $(L - \tfrac{1}{2})^n$. Hence,
	\begin{align*}
		\frac1{|P|}\int_P f(y)\,dy & = 2\frac{|P\cap E|}{|P|}\\
		&\le  2\left(1 - \frac{1}{2L}\right)^n.
	\end{align*}
	This expression is strictly increasing with respect to $L$. For $L \le 1$, its maximum value is attained at $L=1$, yielding $2(1-\tfrac{1}{2})^n = 2^{1-n}$. For $L > 1$, the measure of the intersection is bounded by $|E|=(\tfrac{1}{2})^n$, and then the average is strictly less than $2(\tfrac{1}{2})^n = 2^{1-n}$. Taking the supremum over all cubes $P$ such that $z_0\in P$, we obtain $M f(z_0)\le 2^{1-n}$. Since $f$ is continuous at $z_0$, the maximal function $M f$ is also continuous at $z_0$ (see for instance \cite{MaximalContinuity}). Because $z_0$ belongs to the closure of $Q_0$, the continuity of $M f$ at $z_0$ guarantees that for every $\epsilon > 0$, the set $	\left\{ x \in Q_0 : M f(x) < 2^{1-n} + \epsilon \right\}$ contains the intersection of $Q_0$ with a neighborhood of $z_0$. Since $z_0$ is a vertex of $Q_0$, this intersection has strictly positive Lebesgue measure. This implies that $\essinf_{x\in Q_0} M f(x) \le 2^{1-n}$, and thus $\essinf_{x\in Q_0} M f(x)=2^{1-n}$. 
	
	We can now estimate the $\operatorname{BLO}_w$ seminorm of $M f$,
	\begin{align*}
		\frac1{w(Q_0)}\int_{Q_0} \bigl(M f-\essinf_{Q_0} M f\bigr)\,w
		&\ge \frac1{w(Q_0)}\int_E \bigl(M f-\essinf_{Q_0} M f\bigr)\,w\\
		&= \frac1{w(Q_0)}\int_E \bigl(2-2^{1-n}\bigr)\,w\\
		&= (2-2^{1-n})\frac{w(E)}{w(Q_0)}\\
		&\ge (2-2^{1-n})2^{-n}.
	\end{align*}
	
	This proves that $X\ge c_n$ with $c_n:=(2-2^{1-n})2^{-n}>0$.
	
	Let now $Q$ be an arbitrary cube. We define the auxiliary weight
	\begin{equation*}
		v:= \left[M\left( \frac{w\chi_Q}{w_Q}\right)\right]^\frac{1}{2}.
	\end{equation*}
	It is well known that $v\in A_1$ with $[v]_{A_1}\le c_n$ (see for instance \cite{GCRdF}). Let us define $b:= \log^+ v = \max ( 0, \log v)$, which satisfies $b\in  \operatorname{BMO}$ and $\|b\|_{\operatorname{BMO}}\le \tfrac{3}{2} \|\log v\|_{\operatorname{BMO}} \le \rho_n$ (see \cite{Grafakos}*{Proposition 6.1.4}). Moreover, since $0\le b\le v$ and $Mv\le [v]_{A_1}v<\infty$ almost everywhere, $Mb$ is not identically infinite, so $b$ is admissible in the definition of $X$. Using the definition of $X$, we obtain 
	\begin{equation*}
		\|Mb\|_{\operatorname{BLO}_w } \le \|b\|_{\operatorname{BMO}} X \le \rho_n X, 
	\end{equation*}
	and hence, 
	\begin{equation}\label{eq: proof char}
		\frac{1}{w(Q)}\int_Q \left( Mb(x) - \essinf_{Q} Mb \right) w(x)\, dx \le \rho_n X.
	\end{equation}
	
	We now estimate the essential infimum term. Fix $x\in Q$, and let $P$ be any cube with $x\in P$. By Jensen's inequality, 
	\begin{align*}
		\frac{1}{|P|}\int_P |b(y)|\, dy = & \frac{1}{|P|}\int_P \log^+ v(y)\,  dy\\
		\le & \frac{1}{|P|}\int_P \log(v(y)+1) \, dy\\
		\le & \log \left( \frac{1}{|P|}\int_P v(y) \, dy + 1 \right)\\
		\le & \log \left( Mv(x) + 1 \right).
	\end{align*}
	Consequently, we have $Mb(x)\le \log \left( Mv(x) + 1 \right)$ for every $x\in Q$. Using the $A_1$ condition and that $v\ge 1$ in $Q$, we have 
	\begin{equation*}
		Mb(x) \le \log ( [v]_{A_1} v(x) +1)\le \log([v]_{A_1}+1) + \log v(x) \le \log c_n +\log v(x)
	\end{equation*}
	for almost every $x\in Q$. Therefore,
	\begin{equation*}
		\essinf_{z\in Q} Mb(z) \le \log c_n + \essinf_{z\in Q} \log (v(z)).
	\end{equation*}
	Using the Kolmogorov inequality together with the weak $(1,1)$ boundedness of the maximal operator (see \cite{GrafakosClassic}*{p. 100}), we obtain
	\begin{align*}
		\essinf_{z\in Q} Mb(z) \le & \log c_n  + \essinf_{z\in Q} \log (v(z))\\
		\le & \log c_n + \frac{1}{|Q|}\int_Q \log v \\
		\le & \log c_n + \log \left(\frac{1}{|Q|}\int_Q v\right) \\
		= & \log c_n + \log \left(\frac{1}{|Q|}\int_Q M\left( \frac{w\chi_Q}{w_Q}\right)^\frac{1}{2}\right) \\
		\le & \log c_n + \log \left( 2 \left\| M\left( \frac{w\chi_Q}{w_Q}\right)\right\|_{L^{1,\infty}(\R^n, \frac{1}{|Q|})}^\frac{1}{2}\right) \\
		\le & \log c_n + \log c_n'\\
		= & \log \alpha_n. 
	\end{align*}
	Since the logarithm is increasing, we may assume $\alpha_n\ge e$. Note that
	\begin{equation*}
		\frac{1}{|Q|} \int_Q \log ^{+}\left(\frac{w(x)}{\alpha_n^2 w_Q}\right) w(x) \, dx = \frac{1}{|Q|} \int_{L_Q} \log ^{+}\left(\frac{w(x)}{\alpha_n^2 w_Q}\right) w(x) \, dx,
	\end{equation*}
	where $L_Q = \{ x\in Q : w(x)\ge \alpha_n^2 w_Q\}$. Hence, for a.e. $x\in L_Q$ we have 
	\begin{align*}
		Mb(x)-\essinf_{z\in Q} Mb(z) \ge & Mb(x) - \log \alpha _n \\
		\ge & b(x) - \log \alpha_n\\
		= & \log v(x) - \log \alpha_n \\
		= & \frac{1}{2} \log \left( \frac{1}{\alpha_n^2} M\left( \frac{w\chi_Q}{w_Q}\right)(x) \right) \\
		\ge & \frac{1}{2} \log \left( \frac{1}{\alpha_n^2} \frac{w(x)}{w_Q} \right)\\
		= &\frac{1}{2} \log^+ \left(  \frac{w(x)}{\alpha_n^2 w_Q} \right). 
	\end{align*}

	Combining this with \eqref{eq: proof char}, we obtain
	\begin{align*}
		\frac{1}{2} \frac{1}{|Q|} \int_Q \log ^{+}\left(\frac{w(x)}{\alpha_n^2 w_Q}\right) w(x) \, dx = &  \frac{1}{|Q|} \int_{L_Q} \frac{1}{2} \log ^{+}\left(\frac{w(x)}{\alpha_n^2 w_Q}\right) w(x) \, dx \\
		\le & \frac{1}{|Q|} \int_{L_Q} \left(Mb(x)-\essinf_{z\in Q} Mb(z) \right) w(x) \, dx \\ 
		\le & \frac{1}{|Q|} \int_{Q} \left(Mb(x)-\essinf_{z\in Q} Mb(z) \right) w(x) \, dx \\ 
		\le & \frac{w(Q)}{|Q|} \rho_n X .
	\end{align*}
	
	We conclude that
	\begin{align*}
		\frac{1}{|Q|} \int_Q M(w\chi_Q)(x)\, dx \simeq_n & \frac{1}{|Q|}\int_Q \left( 1+ \log^+ \left( \frac{w(x)}{w_Q} \right) \right) w(x)\, dx \\
		= & \frac{w(Q)}{|Q|} + \frac{1}{|Q|}\int_Q  \log^+ \left( \alpha_n^2 \frac{w(x)}{\alpha_n^2 w_Q} \right) w(x)\, dx \\
		\le & 3\log\alpha_n \frac{w(Q)}{|Q|} + \frac{1}{|Q|}\int_Q  \log^+ \left(  \frac{w(x)}{\alpha_n^2 w_Q} \right) w(x)\, dx \\
		\le & 3\log\alpha_n \frac{w(Q)}{|Q|} +2 \frac{w(Q)}{|Q|} \rho_n X\\
		 \le &  C_n \frac{w(Q)}{|Q|}  X,
	\end{align*}
	where for the first equivalence we refer to \cite{Stein} for the upper bound and to \cite{HP}*{Lemma 6.2} for the lower bound (see also \cite{ClarosRela}*{Remark 3.1 and Lemma 3.4}), and in the last inequality we have used the bound $X\ge c_n>0$ proved above. This completes the proof. 
\end{proof}

We will now apply this characterization of the $A_\infty$ class to prove Corollary \ref{cor: char}, which follows as a direct consequence of the quantitative inequalities established above.

\begin{proof}[Proof of Corollary \ref{cor: char}]
	The first part follows directly from Theorem \ref{thm: 1}. Let us now prove the second part. Let $w$ be a weight and assume that inequality \eqref{eq: cor 2} holds for every nonconstant locally integrable function $f$ such that $Mf$ is not identically infinite and $\|f\|_{\operatorname{BMO}}=1$. In particular, if $f\in \operatorname{BMO}$, then taking $p=1$ in \eqref{eq: cor 2} we have   
	\begin{equation*}
		\|Mf\|_{\operatorname{BLO}_w } \le C \|f\|_{\operatorname{BMO}}.
	\end{equation*}
	By Theorem \ref{thm: characterization Ainfty}, this yields $w\in A_\infty$. In addition, by taking the supremum over $f$ we obtain,
	\begin{align*}
		[w]_{A_\infty} \simeq_n & \sup_{f} \sup_Q \frac{1}{w(Q)}\int_Q \left( Mf(x) - \essinf_{Q} Mf \right) w(x)\, dx \le C.
	\end{align*}
	This completes the proof.
\end{proof}

\section*{Acknowledgment}
Part of this work was carried out while the author was visiting the Department of Mathematics at the Universidad de Buenos Aires. The author wishes to express his sincere gratitude to Ezequiel Rela for his hospitality during this visit.

\section*{Conflicts of Interest}
The author has no conflicts of interest to declare.

\section*{Data Availability Statement}
Data sharing is not applicable to this article, as no datasets were generated or analyzed during the current study.

\end{document}